\documentclass[11pt]{article}
\usepackage{euler}
\usepackage{ams}
\usepackage{xypic}
\newtheorem{dn}{Definition}

\newtheorem{sz}{Proposition}
\newtheorem{tm}{Theorem}

\input xy
\xyoption{all}

\newcommand{\DFT} [0] {{\rm DFT}}

\newcommand{\ADFT} [0] {{\rm ADFT}}

\newcommand{\FADFT} [0] {{\rm FADFT}}
\newcommand{\Sp} [0] {{\rm Tr}}

\begin{document}
\begin{center}{\Large{\bf Slim Normal Bases and Basefield Transforms}}\\
\vspace{.5cm}
{Bj\"orn Grohmann\footnote{email: nn@mhorg.de}}\\
{February 2006}\\
\end{center}
\vspace{.5cm}
\begin{abstract}
In this article, we define the notion of slim (normal) bases and show their existence for various fields. As an application,
an algorithm will be given, that computes the spectrum of a basefield transform by merely using $O(n)$ additions.
\end{abstract}
\section{Introduction}
Let $K\subseteq \C$ be a field, $\zeta_n:=e^{\frac{2\pi i}{n}}$ a primitive $n$-th root of unity and
$G_n:=G(K(\zeta_n)/K)$ the Galois group of the extension $K(\zeta_n)/K$. An element $\vartheta\in K(\zeta_n)$ is called
a {\em normal basis generator} ({\rm NBG}), if the set $\vartheta^{G_n}:=
\{\vartheta^\sigma, \sigma\in G_n\}$ defines a basis of the $K$-vector space $K(\zeta_n)$. In this case,
the set $\vartheta^{G_n}$ is called a {\em normal basis} of the extension $K(\zeta_n)/K$ and  its existence is
guaranteed by the Normal Basis Theorem (cf. \cite{lang}).\\\\
From now on, let $\vartheta$ be a {\rm NBG} of $K(\zeta_n)/K$. Since we are dealing with a separable extension, the space of
functionals on $K(\zeta_n)$ (i.e. ${\rm Hom}_K(K(\zeta_n),K)$) admits a unique basis $\{\phi_{\vartheta^*}^\sigma,\sigma \in G_n\}$, with 
$\phi_{\vartheta^*}^\sigma(\vartheta^\varrho)=\Sp(\vartheta^{*\sigma}\vartheta^\varrho)=\delta_{\sigma,\varrho}$,
for all $\sigma,\varrho\in G_n$. We call the set $\vartheta^{*G_n}:=\{\vartheta^{*\sigma},\sigma\in G_n\}$ the {\em dual basis} of
$\vartheta^{G_n}$. Clearly, this is a normal basis.\\\\
As a member of the family of basefield transforms, the {\em Algebraic Discrete Fourier Transform} (\ADFT), is defined as a linear transfom
$\ADFT_{n,\vartheta} : K^n\longrightarrow K^n,$
with transformation matrix
\begin{equation}
A_{n,\vartheta}:=\left(\phi_{\vartheta^*}(\zeta_n^{kl})\right)_{k,l}=\left(\Sp(\vartheta^*\zeta_n^{kl})\right)_{k,l},
\end{equation}
where $k,l\in\{0,1,\dots,n-1\}$, and it is shown in \cite{beth83}, that this transform computes the spectrum of the
{\em Discrete Fourier Transform} (\DFT) relative to the normal basis $\vartheta^{G_n}$.\\\\
{\bf Example.}\hspace{0.3cm}For $K:=\C$ and $\vartheta:=1$, we get $\ADFT_{n,1}=\DFT_n$.\\\\
{\bf Example.}\hspace{0.3cm}If we take $K:=\R$ and $\vartheta:=(1+i)/2$, the transformation matrix becomes 
$A_{n,\vartheta}=\left(\cos(2\pi kl/n)+\sin(2\pi kl/n)\right)_{k,l}$
and the $\ADFT_{n,\vartheta}$ turns out to be equal to the
{\em Discrete Hartley Transform}.\\\\
{\bf Example.}\hspace{0.3cm}Now take $K:=\Q$, $\vartheta:=(1+i)/2$ and $n:=4$. Then, after some permutation of the columns, the 
transformation matrix $A_{4,\vartheta}$ may be written as
\begin{eqnarray*}
\left(\begin{array}{rrrr}
1&1&1&1\\
1&-1&1&-1\\
1&1&-1&-1\\
1&-1&-1&1
\end{array}
\right)
=\left(\begin{array}{rr}1&1\\1&-1\end{array}\right)\otimes \left(\begin{array}{rr}1&1\\1&-1\end{array}\right),
\end{eqnarray*}
and we see that the $\ADFT_{4,\vartheta}$ can be computed by using merely additions.\\\\
This last example gives rise to the definiton of a slim basis:
\vspace{0.5cm}
\begin{dn}
A normal basis $\vartheta^{G_n}$ will be called a {\em slim basis} of the extension $K(\zeta_n)/K$, if for all
$s\in \Z$ we have
\begin{equation}
\Sp(\vartheta^*\zeta_n^s)\in\{0,\pm 1\}.
\end{equation}
\end{dn}
\vspace{0.5cm}
{\bf Example.}\hspace{0.3cm}It is not hard to see that if $n$ is prime, then the set $\zeta_n^{G_n}$ is a slim
basis of the extension $\Q(\zeta_n)/\Q$, and as we will see in the next section, the same holds if $n$ is {\em squarefree},
i.e. $\mu(n)^2=1$, where $\mu$ denotes the {\em M\"obius function}. Note that this is no longer true when $\mu(n)=0$, since 
in this case the element $\zeta_n$ is not a {\rm NBG} of $\Q(\zeta_n)/\Q$.\\\\
One aim of this paper is to prove the following theorem:
\vspace{0.5cm}
\begin{tm}\label{sbt}
If $G_n\simeq (\Z/n\Z)^\times$, then the extension $K(\zeta_n)/K$ admits a slim basis.
\end{tm}
\vspace{0.5cm}
It is, indeed, enough to prove this in the case $\Q(\zeta_n)/\Q$, since the condition on the Galois group assures that any slim basis of $\Q(\zeta_n)/\Q$
defines as well a slim basis of $K(\zeta_n)/K$.\\\\
The proof of the theorem will proceed in two steps. In the next section, we will reduce to the case of a prime power, and after having done this, an
explicit construction of slim bases will be carried out. The last section concludes the discussion by providing an algorithm, 
that computes the spectrum of the $\ADFT_n$, for $n=2^t$, with the help of $O(n)$ additions in $K$.\\\\
\section{Composition and Decomposition}
{\bf Notation.}\hspace{0.3cm} We keep the notation of the last section. For matrices $A=\left(a_{kl}\right)\in K^{n\times n}$ 
and $B\in K^{m\times m}$, we 
will denote their {\em Kronecker produkt} by $A\otimes B:=\left(a_{kl}B\right)\in K^{nm\times nm}$. In case of $n=m$, we shall
further write $A\sim B$, if there exist permutation matrices $P$ and $Q$, such that $PAQ=B$.\\\\
Now, for positive integers $n_1,n_2$ and $j\in \{1,2\}$, let  $\vartheta_j$ denote a {\rm NBG} of the extension $K(\zeta_{n_j})/K$.
The following proposition shows how to construct bigger transforms by smaller ones:
\vspace{0.5cm}
\begin{sz}\label{composition}
If $(n_1,n_2)=1$ and $K(\zeta_{n_1})\cap K(\zeta_{n_2})=K$, then
\begin{equation}
A_{n_1,\vartheta_1}
\otimes A_{n_2,\vartheta_2}\sim
A_{n_1n_2,\vartheta_1\vartheta_2}.
\end{equation}
\end{sz}
\vspace{0.5cm}
Before we turn to the proof, we state the {\em decomposition property} of a normal basis of the extension $K(\zeta_{n_1},\zeta_{n_2})/K$.
For simplicity, we define $\Sp_j(\cdot):=
\Sp_{K(\zeta_{n_1},\zeta_{n_2})/K(\zeta_{n_j})}(\cdot)$.
\vspace{0.5cm}
\begin{sz}\label{sznbzerleg}
Let $K(\zeta_{n_1})\cap K(\zeta_{n_2})=K$ and
$\vartheta$ be a {\rm NBG} of $K(\zeta_{n_1},\zeta_{n_2})/K$. Then $\vartheta$ admits a decomposition
\begin{equation}
\vartheta=\vartheta_1\vartheta_2,
\end{equation}
where $\vartheta_j\in K(\zeta_{n_j})$, if and only if
\begin{equation}
\Sp(\vartheta)\vartheta=\Sp_1(\vartheta)\Sp_2(\vartheta).
\end{equation}
In particular, if $\Sp(\vartheta)=1=(n_1,n_2)$ and $\vartheta$ can be written as a produkt in the above sense, then
\begin{equation}
A_{n_1n_2,\vartheta}\sim A_{n_1,\Sp_1(\vartheta)}\otimes
A_{n_2,\Sp_2(\vartheta)}.
\end{equation}
\end{sz}
\vspace{0.5cm}
{\bf Proof.} We start with the first proposition. Since $K(\zeta_{n_1})\cap K(\zeta_{n_2})=K$, Galois theory gives
\begin{eqnarray*}
G(K(\zeta_{n_1},\zeta_{n_2})/K) \simeq G(K(\zeta_{n_1})/K)\times
G(K(\zeta_{n_2})/K),
\end{eqnarray*}
and therefore $\vartheta_1\vartheta_2$ is a {\rm NBG} of $K(\zeta_{n_1},\zeta_{n_2})/K$. We further have 
$(\vartheta_1\vartheta_2)^*=\vartheta_1^*\vartheta_2^*$,
by uniqueness of the dual basis. Next, the condition $(n_1,n_2)=1$ leads to $K(\zeta_{n_1},\zeta_{n_2})=K(\zeta_{n_1n_2})$, and
to integers $s_1,s_2$ with $s_1n_1+s_2n_2=1$, so that in summary:
\begin{eqnarray*}
\Sp((\vartheta_1\vartheta_2)^*\zeta_{n_1n_2}) = \Sp_2(\vartheta_1^*\zeta_{n_1}^{s_2})
\Sp_1(\vartheta_2^*\zeta_{n_2}^{s_1})
= \Sp_{K(\zeta_{n_1})/K}(\vartheta_1^*\zeta_{n_1}^{s_2})
\Sp_{K(\zeta_{n_2})/K}(\vartheta_2^*\zeta_{n_2}^{s_1}).
\end{eqnarray*}
For the decomposition, we first note that the element $\Sp_j(\vartheta)$ is easily seen to be a {\rm NBG} of 
the extension $K(\zeta_{n_j})$. On the other hand, if $\vartheta=\vartheta_1\vartheta_2$, then 
$\Sp_1(\vartheta)=\vartheta_1\Sp_1(\vartheta_2)$, and since $K(\zeta_{n_1})\cap K(\zeta_{n_2})=K$, we have $\Sp_1(\vartheta_2)\in K$,
so in particular $\Sp(\vartheta)=\Sp_2(\Sp_1(\vartheta))=
\Sp_1(\vartheta_2)\Sp_2(\vartheta_1)$. We finaly obtain 
$\Sp_1(\vartheta)\Sp_2(\vartheta)=\vartheta_1\vartheta_2
\Sp_1(\vartheta_2)\Sp_2(\vartheta_1)=\vartheta
\Sp(\vartheta)$,
which proves the statement of the second proposition.\hfill{$\Box$}\\\\
{\bf Remark.}\hspace{0.3cm} There are normal bases, that do not decompose in the above sense: the element $\zeta_3\zeta_5+1$ generates
a normal basis of the extension $\Q(\zeta_ {15})/\Q$, but $9(\zeta_3\zeta_5+1)\not=(4-\zeta_3)(2-\zeta_5)$.\\\\
By proposition \ref{composition} we may restrict the search for slim bases to the case where $n$ is a prime power. In addition, by what
has been said before, the ``tame'' case of theorem \ref{sbt} easily follows:
\vspace{0.5cm}
\begin{tm}
Let $n$ be a positive integer with $\mu(n)^2=1$. Then the element $\vartheta:=\mu(n)\zeta_n$ generates a slim basis of the extension
$\Q(\zeta_n)/\Q$.
\end{tm}
\vspace{0.5cm}
\section{Normal bases and Gauss sums}
{\bf Notation.}\hspace{0.3cm} In this section, let $K$ be an abelian number field with Galois group $G:=G(K/\Q)$. 
For a character $\chi$ of $K$, we shall denote its conductor by $f_\chi$ and its {\em field of values} by 
$\Q(\chi):=\Q(\chi(\sigma),\sigma\in G)$. Further, we define for 
$\varrho\in G(\Q(\chi)/\Q)$: $\chi^{\varrho}(\sigma):=(\chi(\sigma))^\varrho$, for all $\sigma\in G$. We will write $\bar\chi$ for the
inverse of the charakter $\chi$, so $\chi\bar\chi=\chi_0$, with $\chi_0(\sigma)=1$ for all $\sigma\in G$.
For an introduction to Dirichlet characters, we refer to \cite{hasse50}, \cite{washington97}.\\\\\
To construct slim bases of the extension $\Q(\zeta_n)/\Q$ for general $n$, we will introduce a representation of algebraic numbers, 
which goes back to H.-W. Leopoldt \cite{leopoldt59}.\\\\
For this, let us recall the definition of a Gauss sum. Let $\chi$ be a Dirichlet character, defined ${\rm mod}\, m$, and $a\in \Z$ be an
integer. The corresponding {\em Gauss sum} is then defined as
\begin{equation}
\tau(\chi|\zeta_{m}^{a}):=
\sum_{s\cap m}\chi(s)\zeta_{m}^{as},
\end{equation}
where ``$s\cap m$'' is short for $s\in (\Z/m\Z)^\times$. If $m=f_\chi$ and $a=1$, we will simply write
\begin{equation}
\tau(\chi):=\tau(\chi|\zeta_{f_{\chi}}).
\end{equation}
Since $\tau(\chi)\tau(\bar\chi)=\chi(-1)f_{\chi}$, this sum is always $\not= 0$. For the general case write 
$\zeta_{m}^{a}=\zeta_{m_{0}}^{a_{0}}$, with $(a_{0},m_{0})=1$, and it holds that $\tau(\chi|\zeta_{m}^{a})=0$, if $f_{\chi}\!\! \not|\, m_{0}$,
and in case of $f_{\chi}|m_{0}$:
\begin{equation}\label{gszer}
\tau(\chi|\zeta_{m}^{a})=\frac{\varphi(m)}{\varphi(m_{0})}
\mu(\frac{m_{0}}{f_{\chi}})\chi(\frac{m_{0}}{f_{\chi}})\bar\chi(a_{0})\tau(\chi),
\end{equation}
where $\varphi$ denotes the {\em Euler phi-function}. For an introduction to Gauss sums and their basic properties, we refer to \cite{hasse50}.\\\\
The following theorem will turn out to be crucial for the construction of slim bases: 
\vspace{0.5cm}
\begin{tm}\label{leo1}
{\rm (Leopoldt)} Let $K$ be an abelian number field with Galois group $G$. For each $\vartheta\in K$ there exists a unique set of
{\em $\chi$-coordinates} $y_K(\chi|\vartheta)$, such that
\begin{equation}\label{dargs}
\vartheta = \frac{1}{|G|}\sum_{\chi}y_{K}(\chi|\vartheta)\tau(\chi),
\end{equation}
where
\begin{equation} \label{wertek}
y_{K}(\chi|\vartheta) \in {\Q}(\chi)
\end{equation}
and
\begin{equation} \label{konj}
y_{K}({\chi}^{\varrho}|\vartheta)=
{y_{K}(\chi|\vartheta)}^{\varrho},
\end{equation}
for all $\varrho\in G(\Q(\chi)/\Q)$.\\\\
On the other hand, any set of coordinates satisfying relations (\ref{wertek}) and (\ref{konj}) gives rise to a unique element of $K$.
\end{tm}
\vspace{0.5cm}
Note that in the original paper there seem to be some minor misprints. In particular the following propostion slightly differs from \cite{leopoldt59}.
\vspace{0.5cm}
\begin{sz}
Let $\vartheta\in K$. Then
\begin{equation}\label{ychif}
y_K(\chi|\vartheta)=\frac{1}{f_\chi}
\sum_\sigma\chi(\sigma)\vartheta^\sigma\,\overline{\tau(\chi)},
\end{equation}
and in particular
\begin{equation}
y_K(\chi|\vartheta^{\sigma})=\bar\chi(\sigma)y_K(\chi|\vartheta).
\end{equation}
\end{sz}
\vspace{0.5cm}
{\bf Proof.} A proof of this proposition (and also of theorem \ref{leo1}) can be found in \cite{mydiss}.\hfill$\Box$\\\\
Now, the next theorem gives the desired connection to normal bases:
\vspace{0.5cm}
\begin{tm}
{\rm (Leopoldt)} The element $\vartheta\in K$ is a {\rm NBG} of the extension $K/\Q$, if and only if $y_K(\chi|\vartheta)\not=0$ for
all characters $\chi$ of $K$.
\end{tm}
\vspace{0.5cm}
\section{Explicit construction of slim bases}
In order to construct slim bases, we start with a calculation of the $\chi$-coordinates of a dual basis: 
\vspace{0.5cm}
\begin{sz}
Let $K$ be an abelian number field and $\vartheta=\frac{1}{|G|}\sum_{\chi}y_K(\chi|\vartheta)\tau(\chi)$ be a {\rm NBG} of
the extension $K/\Q$. Then the $\chi$-coordinates of the dual $\vartheta^*$ are determined by
\begin{equation}\label{spurdual}
y_K(\chi|\vartheta^*)=\frac{\chi(-1)}{f_{\chi}y_K(\bar\chi|\vartheta)}.
\end{equation}
\end{sz}
\vspace{0.5cm}
{\bf Proof.} First, note that the right hand side of (\ref{spurdual}) satisfies the conditions (\ref{wertek}) and (\ref{konj}), since
conjugated characters have the same conductor.
Now, for $\alpha,\beta \in K$ and $\sigma,\varrho\in G$, consider the produkt
\begin{eqnarray*}
\alpha^{\sigma}\beta^{\varrho}=\frac{1}{|G|^2}
\sum_{\chi,\psi}\bar\chi(\sigma)\bar\psi(\varrho)y_K(\chi|\alpha)
y_K(\psi|\beta)\tau(\chi)\tau(\psi).
\end{eqnarray*}
Taking the trace on both sides and rearranging the terms gives:
\begin{eqnarray*}
\Sp(\alpha^{\sigma}\beta^{\varrho})=
\frac{1}{|G|}
\sum_{\chi,\psi}\bar\chi(\sigma)\bar\psi(\varrho)y_K(\chi|\alpha)
y_K(\psi|\beta)\tau(\chi)\tau(\psi)\frac{1}{|G|}\sum_{\pi}
\bar\chi(\pi)\bar\psi(\pi).
\end{eqnarray*}
Next, we notice that the sum $\sum_{\pi}\bar\chi(\pi)\bar\psi(\pi)$ equals $0$, if $\psi\not=\bar\chi$. In case of $\psi=\bar\chi$ it is 
equal to $|G|$, and since $\tau(\chi)\tau(\bar\chi)=\chi(-1)f_{\chi}$,
this leads to
\begin{equation} \label{spf}
\Sp({\alpha}^{\sigma}{\beta}^{\varrho})=
\frac{1}{|G|}\sum_{\chi}y_K(\chi|\alpha)y_K(\bar\chi|\beta)
\bar\chi(\sigma)\chi(\varrho)\chi(-1)f_{\chi}.
\end{equation}
Finally, by taking $\alpha:=\vartheta$ and $\beta:=\vartheta^*$, the proposition follows.\hfill$\Box$\\\\
As an immediate consequence we have:
\vspace{0.5cm}
\begin{sz}
The transformation matrix of the rational $\ADFT$ satisfies
\begin{equation}\label{adftdar}
A_{n,\vartheta}=\frac{1}{\varphi(n)}\left(\sum_{\chi}
\frac{y(\chi|\zeta_{n}^{kl})}{y(\chi|\vartheta)}\right)_{k,l}.
\end{equation}
\end{sz}
\vspace{0.5cm}
{\bf Proof.} This follows from the definition of the matrix $A_{n,\vartheta}$ in combination with (\ref{spurdual}) 
and (\ref{spf}).\hfill$\Box$\\\\
Next, we will compute the $\chi$-coordinates of the roots of unity:
\vspace{0.5cm}
\begin{sz}
For integers $k,l\in \Z$, write $\zeta_n^{kl}=\zeta_{n_0}^{a_0}$, with $(n_0,a_0)=1$. Then
\begin{equation}
y(\chi|\zeta_{n}^{kl})=
y(\chi|\zeta_{n_{0}}^{a_{0}})=\left\{\begin{array}{ll}
\frac{\varphi(n)}{\varphi(n_{0})}\mu(\frac{n_{0}}{f_{\chi}})
\chi(\frac{n_{0}}{f_{\chi}})\bar\chi(a_{0}),& f_\chi|n_0\\
\,0,& {\rm else}\\
\end{array}\right..
\end{equation}
\end{sz}
\vspace{0.5cm}
{\bf Proof.} This follows from (\ref{ychif}) in combination with (\ref{gszer}).\hfill$\Box$\\\\
By what has been said before, we may now assume that $n=p^s$, with $p$ being a prime. We further write
$c=c_{a_0,n_0}=c_{k,l}:=\frac{1}{\varphi(n)}\sum_{\chi}\frac
{y(\chi|\zeta_{n}^{kl})}{y(\chi|\vartheta)}$, and are now led to the following cases:\\\\
$\underline{n_0=1}$: This leads to $c=\frac{1}{y(\chi_0|\vartheta)}$. We set
\begin{equation}
y(\chi_0|\vartheta)=1
\end{equation}
and  obtain $c=1$.\\\\
$\underline{n_0=2}$: Due to the last case: $c=\frac{-1}
{y(\chi_0|\vartheta)}=-1$.\\\\
$\underline{n_0=p\not=2}$: It follows that $c=\frac{1}{\varphi(p)}\left(
-1+\sum_{\chi,f_{\chi}=p}\frac{\bar\chi(a_0)}{y(\chi|\vartheta)}\right)$.
For characters $\chi$, with $f_{\chi}=p$, we set
\begin{equation}
y(\chi|\vartheta)=-1
\end{equation}
and now $c=\frac{-1}{p-1}\sum_{\chi}\chi(a_0)$, so $c=-1$, if
$a_0\equiv 1 \,{\rm mod}\,p$,
else $c=0$.\\\\
$\underline{n_0=2^s,s>1}$: In this case we have $c=\frac{1}{\varphi(2^s)}
\sum_{\chi,f_{\chi}=2^s}\frac{\bar\chi(a_0)}{y(\chi|\vartheta)}$. If we set
\begin{equation}
y(\chi|\vartheta)=\frac{1}{2},
\end{equation}
for characters with conductor $2^s\geq 4$, then, since
$\sum_{\chi,f_{\chi}=2^s}\bar\chi(a_0)=\sum_{\chi}\bar\chi(a_0)-\sum
_{\chi,f_{\chi}<2^s}\bar\chi(a_0)$, we get $c=1$, if
$a_0 \equiv 1 \,{\rm mod}\,2^s$, $c=-1$, if $a_0 \equiv 1 + 2^{s-1}
\,{\rm mod}\,2^{s}$, and else
$c=0$.\\\\
$\underline{n_0=p^s, s > 1, p\not= 2}$: For characters $\chi$ with $f_{\chi}=p^s$,
we set
\begin{equation}
y(\chi|\vartheta)=\frac{1}{p}\sum_{i=1}^{p-1}\left(\frac{i}{p}\right)
\bar\chi(1+ip^{s-1}),
\end{equation}
where
$\left(\frac{i}{p}\right)$ denotes the {\em Legendre symbol}. Since $f_{\chi}=p^{s}$,
$\bar\chi(1+i{p}^{s-1})$ is a primitve $p$-th root of unity
and therefore
$y(\chi|\vartheta)\overline{y(\chi|\vartheta)}=
\frac{1}{p}$. This leads to
\begin{eqnarray*}
c&=&\frac{1}{\varphi(p^{s})}\sum_{i}\left(\frac{i}{p}\right)
\sum_{\chi,f_{\chi}=p^{s}}\bar\chi(a_{0})\chi(1+i{p}^{s-1}).   
\end{eqnarray*}
Now, $\sum_{\chi,f_{\chi}=p^{s}}\bar\chi(a_{0})\chi(1+i{p}^{s-1})\not=0$, if and only if
 $a_{0}\equiv 1+j{p}^{s-1}\,{\rm mod}\,p^s$, $(j,p)=1$, and in this case, we obtain
$c=\left(\frac{j}{p}\right)$.
Else $c=0$.\\\\
As an easy task, it remains to verify, that the defined coordinates $y(\chi|\vartheta)$ fulfill the requirements of
(\ref{wertek}) and (\ref{konj}). This completes the proof of theorem \ref{sbt}.\hfill$\Box$\\\\
\section{A very fast Fourier transform}

As an application we will give an algorithm, that computes the spectrum of the (rational) $\ADFT_n$, for $n=2^t$.
For this, let
\begin{equation}\label{sb2k}
\vartheta_n:=\frac{1}{\varphi(n)}\bigg(1+\frac{1}{2}
\sum_{\chi\not= \chi_0}\tau(\chi)\bigg)
\end{equation}
be the (slim) normal basis generator of the extension $\Q(\zeta_n)/\Q$, that was defined in the last section. Further, let
\begin{equation}
N_{n,\vartheta_n}:=\bigg(\Sp(\vartheta_{n}^{*\sigma}\zeta_{n}^{l})\bigg)_{l,\sigma},
\end{equation}
with $l=0,\dots,n/2-1$ and $\sigma\in G(\Q(\zeta_n)/\Q)$. Since 
$\zeta_n^l=\sum_\sigma\Sp(\vartheta_{n}^{*\sigma}\zeta_{n}^{l})\vartheta_n^\sigma$,
this matrix describes the change between the polynomial
basis $\{\zeta_n^l\}_l$ and the normal basis $\{\vartheta_n^\sigma\}_\sigma$. Writing $A_n:=A_{n,\vartheta_n}$, $N_n:=N_{n,\vartheta_n}$ and
$E_n$ for the {\em identity matrix} of rank n, we have:
\vspace{0.5cm}
\begin{sz}\label{propzer}
A permutation of the columns of $A_n$ leads to the matrix
\begin{equation}\label{saftzer}
\tilde A_n = (A_2 \otimes E_{n/2})\left(
\begin{array}{cc}A_{n/2}&\\&N_{n}\end{array}\right).
\end{equation}
\end{sz}
\vspace{0.5cm}
{\bf Proof.} The permutation of the columns works as follows: first we take the even ones $0, 2, \dots, n-2$, and after that
the inverse of the odd columns (${\rm mod}\, n$) $1, 3^{-1}, \dots, (n-1)^{-1}$. For $n\geq 4$ and any {\rm NBG} $\theta\in\Q(\zeta_n)$, we
are led to a matrix
\begin{eqnarray*}
\tilde A_{n,\theta} = 
\left(
\begin{array}{cc}
A_{n/2,\theta^{\prime}}&\phantom{-}N_{n,\theta}\\
A_{n/2,\theta^{\prime}}&-N_{n,\theta}
\end{array}
\right)=
(A_{2} \otimes E_{n/2})\left(
\begin{array}{cc}A_{n/2,\theta^{\prime}}&\\&N_{n,\theta}\end{array}\right),
\end{eqnarray*}
where the new {\rm NBG} $\theta^{\prime}\in\Q(\zeta_{n/2})$ is derived from $\theta$ by deleting the highest $\chi$-coordinates, i.e
$y(\chi|\theta^{\prime}):=y(\chi|\theta)$, for $f_\chi\leq n/2$. Finaly, setting $\theta:=\vartheta_n$, finishes the proof.\hfill$\Box$\\\\\
A similar proposition can be stated for the matrix $N_n$.
\vspace{0.5cm}
\begin{sz}
Let $n\geq 4$. A permutation of the rows of $N_n$ leads to the marix
\begin{equation}
\tilde N_n = \left(\begin{array}{cc}N_{n/2}&\\&E_{n/4}\end{array}\right)
(A_2\otimes E_{n/4}).
\end{equation}
\end{sz}
\vspace{0.5cm}
{\bf Proof.} The rows are ordered by taking the even ones first, and then the odd ones. Since, for $k\geq 2$ and odd $s$, it
holds that $(2^{k-1}+s)^{-1}\equiv 2^{k-1}+ s^{-1}\,{\rm mod}\, 2^k$, we obtain, for any {\rm NBG} $\theta\in\Q(\zeta_n)$,
\begin{eqnarray*}
\tilde N_{n,\theta} =
\left(
\begin{array}{cc}
N_{n/2,\theta^{\prime}}&\phantom{-}N_{n/2,\theta^{\prime}}\\
R_{n,\theta}&-R_{n,\theta}
\end{array}\right)= 
\left(\begin{array}{cc}N_{n/2,\theta^{\prime}}&\\&R_{n,\theta}\end{array}\right)
(A_2\otimes E_{n/4}),
\end{eqnarray*}
where $R_{n,\theta}$ denotes the matrix
\begin{eqnarray*}
R_{n,\theta}:=\left(\Sp(\theta^*\zeta_n^{ls^{-1}})\right)_{l,s},
\end{eqnarray*}
with $l,s=1,\dots,n/2-1$; $l,s$ odd. Again, the {\rm NBG} $\theta^{\prime}\in\Q(\zeta_{n/2})$ is derived from $\theta$, in the way that
was already discussed in the proof of proposition \ref{propzer}. For $\theta:=\vartheta_n$, this leads to $R_{n,\vartheta_n}=E_{n/4}$, cf.
the last section.\hfill$\Box$\\\\
Putting the pieces together, we obtain:
\vspace{0.5cm}
\begin{tm}
For $n=2^t$, the above defined algorithm computes the spectrum of the $\ADFT_{n,\vartheta_n}$ by merely using $O(n)$ additions in $\Q$.
\end{tm}
\vspace{0.5cm}
{\bf Proof.} First we note that the number of operations needed to perform the base change is linear in $n$:
$\Upsilon_N(n)=\Upsilon_N(n/2)+n/2\leq n-2$.
Therefore, the total number of operations needed by the algorithm is
\begin{eqnarray*}
\Upsilon_A(n)=\Upsilon_A(n/2)+\Upsilon_N(n)+n\in O(n),
\end{eqnarray*}
and since we are working with a slim basis, the theorem follows.\hfill$\Box$\\\\
Figure \ref{sacht} gives an example of the algorithm with $n=8$.
%
%
\begin{figure}
\let\footnotesize=\small\centerline{
\thicklines%
\begin{picture}(350,155)
\put(2,158){\makebox(0,0){$c_0$}}
\put(52,158){\makebox(0,0){$c_1$}}
\put(102,158){\makebox(0,0){$c_2$}}
\put(152,158){\makebox(0,0){$c_3$}}
\put(202,158){\makebox(0,0){$c_4$}}
\put(252,158){\makebox(0,0){$c_5$}}
\put(302,158){\makebox(0,0){$c_6$}}
\put(352,158){\makebox(0,0){$c_7$}}
\multiput(-2,145)(50,0){8}{$\bullet$}
\multiput(0,6)(50,0){8}{\circle{8}}
\multiput(0,0)(50,0){4}{\line(0,1){143}}
\multiput(-4,6)(50,0){8}{\line(1,0){8}}
\multiput(200,0)(50,0){4}{\line(0,1){120}}
\multiput(200,128)(50,0){4}{\circle{16}}
\multiput(194,125)(50,0){4}{{\footnotesize $-1$}}
\multiput(200,136)(50,0){4}{\line(0,1){7}}
\multiput(200,10)(50,0){4}{\line(-3,2){200}}
\multiput(0,10)(50,0){4}{\line(3,2){200}}
\end{picture}%
}%
%
%
\centerline{
\thicklines%
\begin{picture}(350,10)
\multiput(0,0)(50,0){8}{\line(0,1){14}}
\end{picture}%
}%
%
%
%
%
\centerline{
\thicklines%
\begin{picture}(350,10)
\multiput(0,0)(50,0){8}{\line(0,1){14}}
\end{picture}%
}%
%
%
\centerline{
\thicklines%
\begin{picture}(200,67)%
\multiput(0,6)(50,0){4}{\circle{8}}
\multiput(0,0)(50,0){2}{\line(0,1){70}}
\multiput(-4,6)(50,0){4}{\line(1,0){8}}
\multiput(100,0)(50,0){2}{\line(0,1){46}}
\multiput(100,54)(50,0){2}{\circle{16}}
\multiput(94,51)(50,0){2}{{\footnotesize $-1$}}
\multiput(100,62)(50,0){2}{\line(0,1){8}}
\multiput(100,10)(50,0){2}{\line(-5,3){100}}
\multiput(0,10)(50,0){2}{\line(5,3){100}}
\end{picture}%
\begin{picture}(150,67)%
\multiput(0,6)(100,0){2}{\circle{8}}
\multiput(0,0)(50,0){2}{\line(0,1){70}}
\multiput(150,0)(50,0){1}{\line(0,1){70}}
%
\multiput(-4,6)(100,0){2}{\line(1,0){8}}
\multiput(100,0)(50,0){1}{\line(0,1){46}}
\multiput(100,54)(50,0){1}{\circle{16}}
\multiput(94,51)(50,0){1}{{\footnotesize $-1$}}
\multiput(100,62)(50,0){1}{\line(0,1){8}}
\multiput(100,10)(50,0){1}{\line(-5,3){100}}
\multiput(0,10)(50,0){1}{\line(5,3){100}}
\end{picture}%
}%
%
%
\centerline{
\thicklines%
\begin{picture}(350,7)
\multiput(0,0)(50,0){8}{\line(0,1){14}}
\end{picture}%
}%
%
%
%
%
\centerline{
\thicklines%
\begin{picture}(350,7)
\multiput(0,0)(50,0){8}{\line(0,1){14}}
\end{picture}%
}%
%
%
\centerline{
\thicklines%
\begin{picture}(100,42)%
\multiput(0,6)(50,0){2}{\circle{8}}
\put(0,0){\line(0,1){43}}
\multiput(-4,6)(50,0){2}{\line(1,0){8}}
\put(50,0){\line(0,1){20}}
\put(50,28){\circle{16}}
\put(44,25){{\footnotesize $-1$}}
\put(50,36){\line(0,1){7}}
\put(50,10){\line(-3,2){50}}
\put(0,10){\line(3,2){50}}
\end{picture}%
\begin{picture}(100,42)%
\multiput(0,6)(50,0){2}{\circle{8}}
\put(0,0){\line(0,1){43}}
\multiput(-4,6)(50,0){2}{\line(1,0){8}}
\put(50,0){\line(0,1){20}}
\put(50,28){\circle{16}}
\put(44,25){{\footnotesize $-1$}}
\put(50,36){\line(0,1){7}}
\put(50,10){\line(-3,2){50}}
\put(0,10){\line(3,2){50}}
\end{picture}%
\begin{picture}(100,42)%
\multiput(0,6)(50,0){2}{\circle{8}}
\put(0,0){\line(0,1){43}}
\multiput(-4,6)(50,0){2}{\line(1,0){8}}
\put(50,0){\line(0,1){20}}
\put(50,28){\circle{16}}
\put(44,25){{\footnotesize $-1$}}
\put(50,36){\line(0,1){7}}
\put(50,10){\line(-3,2){50}}
\put(0,10){\line(3,2){50}}
\end{picture}%
\begin{picture}(50,42)%
\multiput(0,6)(50,0){2}{\circle{8}}
\put(0,0){\line(0,1){43}}
\multiput(-4,6)(50,0){2}{\line(1,0){8}}
\put(50,0){\line(0,1){20}}
\put(50,28){\circle{16}}
\put(44,25){{\footnotesize $-1$}}
\put(50,36){\line(0,1){7}}
\put(50,10){\line(-3,2){50}}
\put(0,10){\line(3,2){50}}
\end{picture}%
}%
%
%
\centerline{
\thicklines%
\begin{picture}(350,60)%
\multiput(0,60)(50,0){8}{\line(0,1){2}}
\multiput(0,60)(100,0){2}{\vector(0,-1){46}}
\multiput(250,60)(100,0){2}{\vector(0,-1){46}}
\multiput(50,60)(100,0){2}{\vector(3,-1){146}}
\multiput(200,60)(100,0){2}{\vector(-3,-1){146}}
\multiput(-3,5)(50,0){8}{$\bullet$}
\put(2,-4){\makebox(0,0){$\hat{c}_0$}}
\put(52,-4){\makebox(0,0){$\hat{c}_1$}}
\put(102,-4){\makebox(0,0){$\hat{c}_2$}}
\put(152,-4){\makebox(0,0){$\hat{c}_3$}}
\put(202,-4){\makebox(0,0){$\hat{c}_4$}}
\put(252,-4){\makebox(0,0){$\hat{c}_5$}}
\put(302,-4){\makebox(0,0){$\hat{c}_6$}}
\put(352,-4){\makebox(0,0){$\hat{c}_7$}}
\end{picture}%
}%
\caption{\label{sacht} $\FADFT_8$}
\end{figure}


\begin{thebibliography}{99}
\bibitem{beth83}
Th. Beth, W. Fumy, R. M\"uhlfeld, {\it Zur Algebraischen Diskreten
Fourier-Transformation}, Arch. Math. 40, 238--244, 1983
\bibitem{mydiss}
B. Grohmann, {\it Zur Arithmetik in abelschen Zahlk\"orpern}, Dissertation, Universit\"at Karlsruhe, 2005,
avail. at http://www.ubka.uni-karlsruhe.de/cgi-bin/psview?document=2005/informatik/25

\bibitem{hasse50}
H. Hasse, {\it Vorlesung \"uber Zahlentheorie}, Die Grundl. der Math.
Wissensch. in Einzeld., Band LIX, Springer Berlin G\"ottingen Heidelberg,
1950

\bibitem{lang}
S. Lang, {\it Algebra}, 3rd ed., Addison-Wesley, 1994

\bibitem{leopoldt59}
H.-W. Leopoldt, {\it \"Uber die Hauptordnung der ganzen Elemente eines
abelschen Zahlk\"orpers}, J. reine u. angew. Math. 201, 119--149, 1959

\bibitem{washington97}
L. C. Washington, {\it Introduction to Cyclotomic Fields}, GTM 83, 2nd Ed.,
Springer New York, 1997

\end{thebibliography}
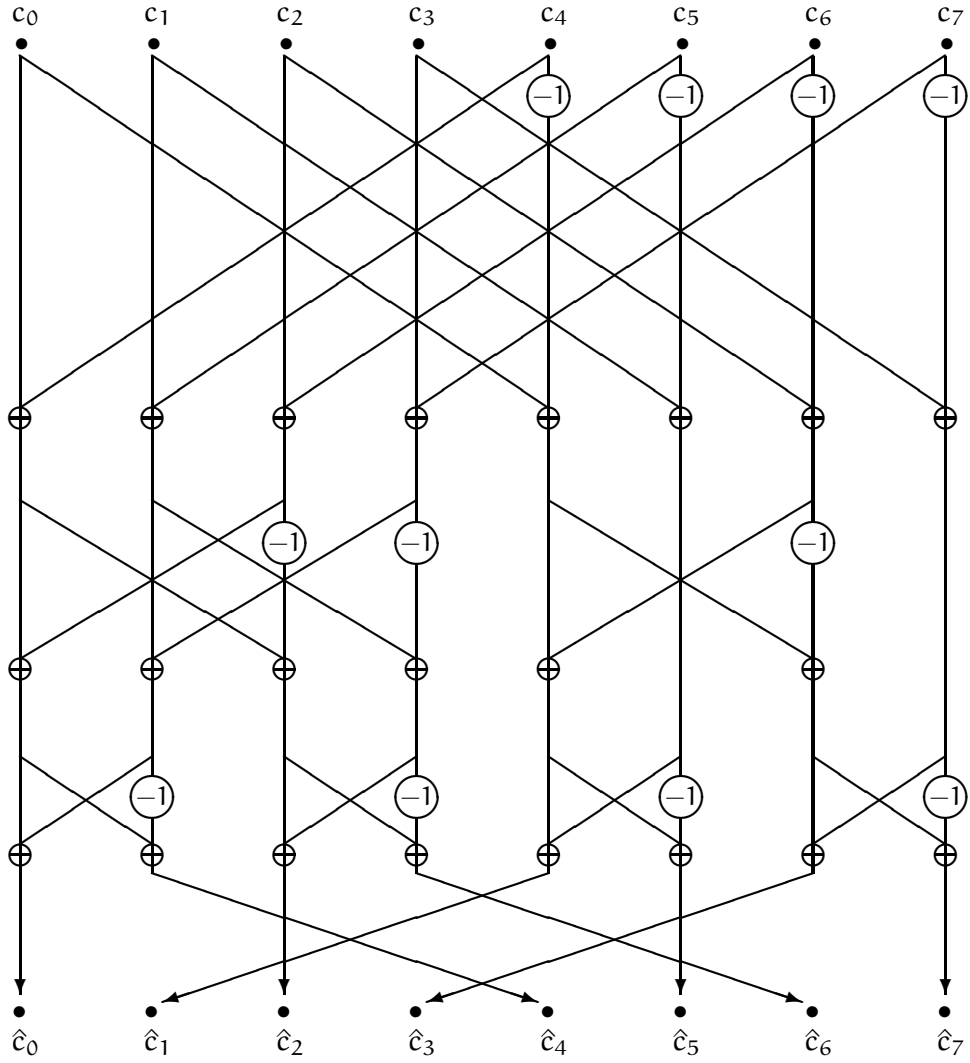
\end{document}